\documentclass[12pt]{amsart}
\usepackage{amssymb}
\usepackage{amsmath}
\newcommand\C{\mathbb C}
\newcommand\R{\mathbb R}
\newcommand\g{{\mathfrak g}}
\newcommand\gl{{\mathfrak{gl}}}

\newcommand\Z{{\mathbb Z}}
\renewcommand\k{{\mathfrak k}}
\renewcommand\Im{\operatorname{Im}}
\newtheorem{Thm}{Theorem}
\theoremstyle{definition}

\title{The Kempf-Ness theorem and Invariant Theory}
\author{Ivan V. Losev}
\begin{document}
\begin{abstract}
We give new proofs of some well-known results from Invariant Theorey
using the Kempf-Ness theorem.
\end{abstract}
\maketitle

\section{Introduction}
This article does not contain any new results. Its goal is to deduce
some well-known results of Invariant Theory from the Kempf-Ness
theorem.

In the sequel $G$ denotes a complex reductive algebraic group. By
a small fraktur letter we denote the Lie algebra of a Lie group
denoted by the corresponding capital Latin letter.

Let us state the results we want to prove.
\begin{Thm}[The Matsushima criterion]\label{Thm:1}
Let $H$ be an algebraic subgroup of $G$. The homogeneous space $G/H$
is affine iff $H$ is reductive.
\end{Thm}

The theorem was proved independently by  Matsushima,
\cite{Matsushima}, and Onishchik~\cite{Onishchik}. It has many
different proofs, see~\cite{Arzhantsev} for references and one more
proof. One part of the theorem (the "if" part) is easy. We give a
new proof of the other part.

\begin{Thm}[The Luna criterion]\label{Thm:2}
Let $H$ be a reductive subgroup of $G$, $X$ an affine $G$-variety
and $x$ an $H$-stable point in $X$. Then $N_G(H)x$ is closed in $X$
iff   $Gx$ is.
\end{Thm}

 The Luna criterion for orbit's closedness
was originally proved in~\cite{Luna2}. The proof is quite involved.
An alternative (and easier) proof was obtained by Kempf,
\cite{Kempf}. Again, the "if" part of the theorem is easy, and we
give a new proof of the difficult part.

\begin{Thm}\label{Thm:3}
Let $H_1,H_2$ be reductive subgroups of $G$. Then the action
$H_1:G/H_2$ is stable, i.e. an orbit in general position is closed.
\end{Thm}

This result also has different proofs, see~\cite{Vinberg} for
details.

The proofs of all three theorems are based on the Kempf-Ness
criterion for the orbit closedness which we state now.

Let $V$ be a $G$-module and $K$ a compact form of $G$. Choose a
$K$-invariant hermitian scalar product $(\cdot,\cdot)$ on $V$.
Define a map $\mu:V\rightarrow \k^*$ by the formula
\begin{equation}\label{eq:0} \langle \mu(v),\xi\rangle=\frac{1}{2i}(
\xi v,v)
\end{equation}

$\mu(v)$ lies in $\k^*$ because the image of $\k$ in $\gl(V)$
consists of skew-hermitian operators. The map $\mu$ is the {\it
moment map} for the action $K:V$.

\begin{Thm}[The Kempf-Ness criterion]\label{Thm:4}
For $v\in V$ the orbit $Gv$ is closed iff $Gv\cap \mu^{-1}(0)\neq
\varnothing$.
\end{Thm}

A stronger result was proved in~\cite{Kempf_Ness}.

The author wishes to thank 
I.V. Arzhantsev for some
remarks on an earlier version of this article.

\section{Proofs}

At first, for convenience of a reader we give

\begin{proof}[Proof of the Kempf-Ness criterion]
Let us note that, again, one implication in the theorem is easy.
Denote by $||v||$ the length of a vector $v\in V$ with respect to
$(\cdot,\cdot)$.  If an orbit $Gv$ is closed, then there is a point
$v_0\in Gv$ such that $||v_0||=\min_{u\in Gv}||u||$. Thus for any
$\xi\in \k$ there is the equality
$\frac{d}{dt}(\exp(it\xi)v,\exp(it\xi)v)|_{t=0}=0$, or,
equivalently, $0=(i\xi v,v)+(v,i\xi v)=2i(\xi v,v)=-4\langle
\mu(v),\xi\rangle $.

Conversely, let $v\in V$ be such that $\mu(v)=0$. Assume that the
orbit $Gv$ is not closed. It follows from the Hilbert-Mumford
criterion that the exist a one-parameter subgroup
$\tau:\C^\times\rightarrow G$ such that $\lim_{t\rightarrow
0}\tau(t)v$ exists and is not equal to $v$ and $\tau$ is compatible
with $K$, i.e. $\tau(z)\in K$ if $|z|=1$. The last statement follows
easily from a proof of the criterion given, for example,
in~\cite{Kraft}.

Let $v=\sum_{i\in \Z} v_i$ be the weight decomposition with respect
to $\tau$. It can be easily seen that $(v_i,v_j)=0$ provided $i\neq
j$. The limit $\lim_{t\rightarrow 0}\tau(t)v$ exists iff $v_i=0$ for
all $i<0$ and is equal to $v_0$. The equality
$\langle\mu(v),\xi\rangle=0$ for $\xi=\frac{d}{dt}\tau|_{t=0}$ can
be rewritten as $\sum_{i\in \Z} i(v_i,v_i)=0$. Hence $v_i=0$ for
$i>0$ and $v=v_0=\lim_{t\rightarrow 0}\tau(t)v$.
\end{proof}

\renewcommand{\proofname}{Proof of Theorem~\ref{Thm:1}}
\begin{proof}
First suppose that $H$ is a reductive subgroup of $G$. Then $G/H$ is
a categorical quotient in the sense of Geometric Invariant Theory,
see, for example,~\cite{VP}. In particular, $G/H$ is affine.

Suppose now that $G/H$ is affine. Then there are a $G$-module $V$
and a closed $G$-equivariant embedding $G/H\hookrightarrow V$ (see
\cite{VP}, $\S$1). Choose a compact form $K\subset G$ and a
$K$-invariant hermitian scalar product $(\cdot,\cdot)$ on $V$.   By
(the easy part of) Theorem~\ref{Thm:4}, one can find a point $v\in
G/H\cap \mu^{-1}(0)$.

The real 2-form $\omega(u,v)=\Im \langle u,v\rangle$ on $V$ is
symplectic. Moreover, for any complex submanifold $X\subset V$ the
restriction of $\omega$ to $X$ is again symplectic. Clearly,
$\omega$ is $K$-invariant. In particular,

\begin{equation}\label{eq:4}
\omega(\xi v_1,v_2)+\omega(v_1,\xi v_2)=0, \forall \xi\in \k, u,v\in
V
\end{equation}

The equality (\ref{eq:4}) applied for $v_1=v, v_2=\eta v,\eta\in
\k,$ implies
\begin{equation}\label{eq:2}
\langle\mu(v), [\xi,\eta]\rangle=\omega(\xi v,\eta v)
\end{equation}

 By (\ref{eq:2}), the orbit $Kv$ is an
isotropic submanifold of $Gv$. Therefore $2\dim_{\R}Kv\leqslant
\dim_\R Gv$ or, equivalently,
\begin{equation}\label{eq:1}\dim_{\R}Kv\leqslant \dim_\C Gv.\end{equation} But
$\dim_\R Kv=\dim_\R K-\dim_\R\k_v=\dim_\C G-\dim_\R (\g_v\cap\k)$.
Thus (\ref{eq:1}) implies $\dim_\R(\g_v\cap \k)\geqslant
\dim_\C\g_v$. It follows that $\g_v\cap\k$ is a real form of $\g_v$.
Therefore $\g_v$ is reductive.
%
\end{proof}

\renewcommand{\proofname}{Proof of Theorem~\ref{Thm:2}}
\begin{proof}
For convenience of a reader we give a proof of an easy part of the
theorem due to Luna. That is, we prove that if $Gx$ is closed, then
 $N_G(H)x$ is. Let $y\in Gx$. We have $T_y((Gx))^H=
(\g/\g_y)^H=\g^H/(\g_y)^H=T_y(G^Hy)=T_y(N_G(H)y)$. It follows that
$(Gx)^H$ is a smooth variety and its components are
$N_G(H)^\circ$-orbits. In particular, for any $y\in X^H$ the orbit
$N_G(H)y$ is closed.

Now we prove that if $N_G(H)x$ is closed, then so is $Gx$. Let us
embed $X$ into a $G$-module $V$.

 One can choose compact forms $K,K_1,K_2$ of $G,$ $N_G(H),H$,
respectively, such that $K_2\subset K_1\subset K$. We choose a
$K$-invariant hermitian scalar product $(\cdot,\cdot)$ on $V$ and
define the moment map $\mu$ by (\ref{eq:0}).

Let $v\in V^H$ be such that $N_G(H)v$ is closed. Denote by $\mu_1$
the moment map for the action $K_1:V^H$ (in~(\ref{eq:0}) we take the
restriction of $(\cdot,\cdot)$ to $V^H$ for the  scalar product).
Choose an invariant  scalar product $\langle\cdot,\cdot\rangle$ on
$\k$ and identify $\k,\k_1$ with their duals via
$\langle\cdot,\cdot\rangle$. Let us prove that $\mu_1=\mu|_{V^H}$.
It can be seen directly, that the map $\mu:V\rightarrow \k$ is
$K$-equivariant. It follows that $\mu(V^H)\subset \k^{K_2}\subset
\k_1$. The equality $\mu_1=\mu|_{V^H}$ follows now directly from the
definitions of $\mu,\mu_1$.

By Theorem~\ref{Thm:4} applied to the action $N_G(H):V^H$, one can
choose a vector $v_1\in N_G(H)v\cap \mu_1^{-1}(0)$. Since
$\mu_1=\mu|_{V^H}$, $v_1\in \mu^{-1}(0)$. Applying
Theorem~\ref{Thm:4} again, we see that the orbit $Gv$ is closed.
\end{proof}

\begin{proof}[Proof of Theorem~\ref{Thm:3}]
Embed $G/H_2$ into a $G$-module $V$, fix a compact form $K\subset
G$,
 a $K$-invariant hermitian scalar product $(\cdot,\cdot)$ on $V$
 and a compact form $K_1$ of $H$ such that $K_1\subset K$. Let $\pi$
 denote the natural projection $\k^*\rightarrow \k_1^*$. The map $\pi\circ\mu$
 is the moment map for the action $K_1:V$. There exists $v\in G/H_2$
 such that $\mu(v)=0$. Note that $\mu(kv)=0$ for all $k\in K$.
In particular, $\pi\circ\mu(Kv)=0$. By Theorem~\ref{Thm:4}, the
orbit $H_1kv$ is closed for any $k\in K$.
 It remains to check that the
 subset $Kv$ is dense in $Gv$. Assume the converse: there exists a
 proper
 closed subvariety $Y\subset G/H_2$ containing $Kv$. Replacing $Y$
 by $\bigcap_{k\in K}kY$ we may assume that $Y$ is $K$-invariant.
 Since $K$ is Zariski-dense in $G$, $Y$ is $G$-invariant.
 Contradiction.
\end{proof}


\begin{thebibliography}{99}
\bibitem[A]{Arzhantsev} I.V. Arzhanysev. {\it Invariant
ideals and the Matsushima criterion}. Preprint (2005),
arXiv:math/AG.0506430.

\bibitem[Ke]{Kempf} G. Kempf. {\it Instability in invariant theory}. Ann.
Math. II. Ser. 108(1978), p. 299-316.

\bibitem[KN]{Kempf_Ness} G. Kempf, L. Ness. {\it The length of vectors in
representation spaces}. Lect. Notes. Math. 732. Springer Verlag,
1979. p. 233-243.

\bibitem[Kr]{Kraft}  H. Kraft. {\it Geometrishe Methoden in der
Invarianttheorie}. Braunschweig/Wiesbaden, Viewveg, 1985.

\bibitem[L]{Luna2} D. Luna. {\it Adh\'{e}rences
d'orbite et invariants}. Invent. Math, 29(1975), p. 231-238.


\bibitem[M]{Matsushima} Y. Matsushima. {\it Espaces
homog\`{e}nes de Stein des groupes des Lie complexes}. Nagoya Math.
J. 16(1960), p. 205-216.

\bibitem[O]{Onishchik} A.L. Onishchik. {\it Complex hulls of complex homogeneous
spaces}. Dokl. Akad. Nauk SSSR, 130(1960), 4, p. 88-91. English
translation: Sov. Math. Dokl., 1(1960), p. 88-91.


\bibitem[PV]{VP} V.L. Popov, E.B. Vinberg. {\it Invariant theory}. Algebraic
geometry IV, Encyclopaedia of Math. Sciences, vol.55. Springer
Verlag, 1994, pp. 123--278.

\bibitem[V]{Vinberg} E.B. Vinberg. {\it On stability of actions of reductive
algebraic groups}. in "Lie algebras, rings and related topics", Fong
Yuen, A.A. Mikhalev, E. Zelmanov eds. Springer-Verlag, Hong Kong
(2000), 188-202.
\end{thebibliography}
\end{document}